\documentclass[11pt,a4paper]{article}
\usepackage{amsmath}
\usepackage{amssymb}
\usepackage{latexsym}

\newcommand{\dual}{\makebox[0mm]{}^{{\scriptstyle\vee}}}
\newtheorem{theorem}{Theorem}[section]

\newtheorem{lemma}[theorem]{Lemma}
\newtheorem{proposition}[theorem]{Proposition}
\newtheorem{definition}[theorem]{Definition}
\newtheorem{corollary}[theorem]{Corollary}

\newtheorem{exmp}[theorem]{Example}
\newtheorem{exmps}[theorem]{Examples}
\newtheorem{rem}[theorem]{Remark}

\newenvironment{remark}{\begin{rem}\rm}{\end{rem}\rm}
\newcommand{\prf}{{\em Proof}. }
\newcommand{\qed}{\hspace*{\fill}$\Box$}

\newcommand{\beeq}[1]{\begin{eqnarray}\label{#1}}
\newcommand{\eneq}{\end{eqnarray}}

\newcommand{\ka}{{\cal A}}
\newcommand{\kb}{{\cal B}}
\newcommand{\kc}{{\cal C}}

\newcommand{\kk}{{\cal K}}

\newcommand{\ko}{{\cal O}}

\newcommand{\kt}{{\cal T}}
\newcommand{\kx}{{\cal X}}

\newcommand{\kz}{{\cal Z}}

\newcommand{\IN}{{\mathbb N}}

\newcommand{\IQ}{{\mathbb Q}}
\newcommand{\IR}{{\mathbb R}}
\newcommand{\IZ}{{\mathbb Z}}

\newcommand{\Pic}{{\rm Pic}}

\newcommand{\id}{{\rm id}}
\newcommand{\im}{{\rm im}}

\newcommand{\codim}{{\rm codim}}

\newcommand{\verylongarrow}[1]{\hbox to #1{\rightarrowfill}}

\begin{document}
{\centerline{\Large\bf The K\"ahler cone of a compact hyperk\"ahler manifold}}

\bigskip

{\centerline{\large\bf Daniel Huybrechts}

\bigskip

\bigskip

\begin{abstract}
This paper is a sequel to \cite{Huybrechts2}. We study a number
of questions only touched upon in \cite{Huybrechts2} in more detail.
In particular: What is the relation between two birational
compact hyperk\"ahler manifolds? What is the shape of the cone
of all K\"ahler classes on such a manifold? How can
the birational K\"ahler cone be described?
Most of the results are motivated by either the well-established
two-dimensional theory, i.e.\ the theory of K3 surfaces, or
the theory of Calabi-Yau threefolds and string theory.
\end{abstract}

\bigskip

\bigskip

\section{Introduction}\label{intro}

Let $X$ be a K3 surface and let $\alpha\in H^{1,1}(X,\IR)$
be a class in the positive
cone with $\alpha.C\ne0$ for any smooth rational curve $C$.
Then there exist smooth rational curves $C_i$ such that for the effective cycle
$\Gamma:=\Delta+\sum C_i\times C_i\subset X\times X$ (where $\Delta$ is
the diagonal) the class $[\Gamma]_*(\alpha)$ is a K\"ahler class.
This is a consequence of the transitivity of the action of the Weyl-group
on the set of chambers in the positive cone and the description of the
K\"ahler cone of a K3 surface due to Todorov and Siu (cf.\ \cite{Periodes}).
In higher dimension the following result was proved in
\cite[Cor.\ 5.2]{Huybrechts2}:

\bigskip
{\it --- Let $X$ be a compact hyperk\"ahler manifold of dimension $2n$
and let $\alpha\in H^{1,1}(X,\IR)$ be
a very general class in the positive cone.
Then there exists another compact hyperk\"ahler manifold $X'$  and
an effective cycle $\Gamma:=Z+\sum Y_i\subset X\times X'$ of pure
dimension $2n$ such that
$X\leftarrow Z\to X'$ defines a birational correspondence,
the projections $Y_i\to X$, $Y_i\to X'$ are of positive dimensions,
and $[\Gamma]_*(\alpha)$ is a K\"ahler class on
$X'$.}

\bigskip

In this article we present a few applications of this result.
In Sect.\ \ref{biratman} we will show the following theorem:

\bigskip

{\it --- Let $X$ and $X'$ be compact hyperk\"ahler manifolds. If $X$
and $X'$ are birational then there exist two smooth proper morphisms
${\cal X}\to S$ and ${\cal X}'\to S$ with $S$ smooth and one-dimensional
such that for a distinguished point $0\in S$ one has
${\cal X}_0\cong X$ and  ${\cal X}'_0\cong X'$ and such that there exists an
isomorphism
${\cal X}|_{S\setminus\{0\}}\cong{\cal X'}|_{S\setminus\{0\}}$.}

\bigskip

In particular, $X$ and $X'$ are deformation equivalent and, hence,
diffeomorphic. This immediately proves that the Betti numbers of $X$
and $X'$ coincide and, more precisely, their Hodge structures are
isomorphic. The result was proved in \cite{Huybrechts2}
for projective hyperk\"ahler manifolds.
The proof given here is different even in the projective case.
The fact about the Betti numbers also follows from a theorem of
Batyrev \cite{Ba} and Kontsevich \cite{DL}. They show that if $X$ and $X'$
are birational smooth projective varieties with trivial canonical bundle,
then $b_i(X)=b_i(X')$ (cf.\ \cite{Ba}) and their Hodge structures
are isomorphic (cf.\ \cite{DL}). For hyperk\"ahler manifolds
the above result proves that the cohomology rings of birational
hyperk\"ahler manifolds are isomorphic. This is no longer true
for arbitrary Calabi-Yau manifolds.

In Sect.\ \ref{Kaehlercone} we give a description of the K\"ahler cone of a
compact hyperk\"ahler manifold very much in the spirit of the known
one for K3 surfaces (Prop.\ \ref{Kcone}):

\bigskip

--- {\it  Let $X$ be a compact hyperk\"ahler manifold.
A class $\alpha\in H^{1,1}(X,\IR)$ is in the closure
$\overline \kk_X$ of the cone of all K\"ahler classes if and only if
$\alpha$ is in the closure $\overline\kc_X$ of the positive cone
and $\alpha.C\geq0$ for any rational curve $C\subset X$.}

\bigskip

The arguments applied to the K3 surface situation differ
from the known ones in as much as they make no use of the
Global Torelli Theorem.

Any irreducible curve of negative self-intersection on a K3 surface
is a smooth rational curve. In higher dimensions one can prove
the following generalization (Prop.\ \ref{Uni}):

\bigskip

--- {\it Let $C$ be an irreducible curve in a compact hyperk\"ahler
manifold $X$. If $q_X([C])<0$, then $C$ is contained
in a uniruled subvariety $Y\subset X$.}

\bigskip

Here, $q_X$ is the quadratic form on $H^{1,1}(X)\cong H^{2n-1,2n-1}(X)$
introduced by Beauville and Bogomolov. Stronger versions of the result are
conjectured in Sect.\ \ref{negcurve}.

\bigskip

Due to the existence of birational maps between compact hyperk\"ahler
manifolds which do not extend to isomorphisms, one has in higher
dimensions the notion of the birational K\"ahler cone.
The birational K\"ahler cone of Calabi-Yau threefolds has been intensively
studied, partly motivated by mirror symmetry, in \cite{Kawamata2, Morrison}.
Analogously to the description of the K\"ahler cone, but replacing
curves by divisors, one has (Prop.\ \ref{BirKcone}):

\bigskip
--- {\it The closure of the birational K\"ahler cone is the set
of all classes $\alpha\in\overline\kc_X$ such that $q_X(\alpha,[D])\geq0$
for all uniruled divisors $D$.}

\bigskip

The present article is the third version of a paper that
had originally been prepared
for publication in the proceedings of the conference on the occasion
of F.\ Hirzebruch's 70th birthday in Warsaw 1998. As all results of
this paper depend on the projectivity criterion Thm.\ 3.11 in
\cite{Huybrechts2} whose proof in \cite{Huybrechts2} is not correct,
it has never been published. In the meantime the projectivity criterion
for hyperk\"ahler manifolds could be proved by using a new
theorem of Demailly and Paun (see the Erratum in \cite{Huybrechts}).
So, all results of this paper are fully proved now. The second
version of the article contained a number of basic facts on the
Beauville-Bogomolov quadratic form. They are of independent interest, but
not directly relevant to the main results of the article.
So we have decided to omit them in this version.

\bigskip

\noindent
{\bf Notations.} In this paper a compact hyperk\"ahler manifold
is a simply connected K\"ahler manifold $X$ such that
$H^0(X,\Omega_X^2)$ is spanned by an everywhere non-degenerate
holomorphic two-form $\sigma$. Equivalently, it is a compact complex
manifold that admits a K\"ahler metric with holonomy ${\rm Sp}(n)$,
where $\dim(X)=2n$. If $X$ is a compact hyperk\"ahler manifold then
there exists a natural quadratic form $q_X$
of index $(3,b_2(X)-3)$ on $H^2(X,\IZ)$ with the property that
$q_X(\alpha)^n=\int_X\alpha^{2n}$ up to a scalar factor (cf.\
\cite{Beauville1},\cite{Bo2}). The {\it positive cone} $\kc_X$
is the connected component of the cone
$\{\alpha\in H^{1,1}(X,\IR)|q_X(\alpha)>0\}$
that contains the {\it K\"ahler cone} $\kk_X$ of all K\"ahler classes.\\
The word ``birational'' is used even when the manifolds in question
are not projective; what is meant in this case, of course,
is ``bimeromorphic''. I hope this will cause no confusion.

\bigskip

\noindent
{\bf Acknowledgement:} I am grateful to Y.\ Kawamata for 
interesting discussions on the content of this paper.
In particular, he pointed out an oversimplification
in the proof of Prop.\ \ref{Fujikimod} in an earlier version.

\bigskip


\section{Birational hyperk\"ahler manifolds}\label{biratman}

Let $f:X' - - \to X$ be a birational map of compact complex
manifolds that induces an isomorphism on open sets whose
complements are analytic subsets of codimension $\geq 2$.
Then $f^*: H^2(X,\IR)\cong H^2(X',\IR)$ and $f^*: \Pic(X)\cong\Pic(X')$,
where $(f^*)^{-1}=(f^{-1})^*$.

If $X$ and $X'$ are projective and the pull-back $f^*L$ of an ample
line bundle $L$ on $X$ is an ample line bundle on $X'$, then,
obviously, $f$ extends to an isomorphism $X\cong X'$. The analytic
analogue was proved by Fujiki \cite{Fujiki3}: If $X$ and $X'$ are K\"ahler
and the pull-back $f^*(\alpha)\in H^2(X',\IR)$ of a K\"ahler
class $\alpha\in H^2(X,\IR)$ is a K\"ahler class on $X'$, then
$f$ extends to an isomorphism $X'\cong X$. A slight modification
of Fujiki's arguments also proves:

\begin{proposition}\label{Fujikimod}---
Let $f:X' - - \to X$ be a birational map of compact complex
manifolds with nef canonical bundles $K_X$ and $K_{X'}$, respectively.
If $\alpha\in H^2(X,\IR)$ is a class such that $\int_C\alpha>0$ and
$\int_{C'}f^*(\alpha)>0$ for all rational
curves $C\subset X$ and $C'\subset X'$, then $f$ extends to an isomorphism
$X\cong X'$.

\end{proposition}

\prf Let $\pi:Z\to X$ be a sequence of blow-ups resolving $f$, i.e.\
there exists a commutative diagram of the form
$$\begin{array}{cccccc}
&&Z&&\\
&\pi'\swarrow&&\searrow\pi&\\
&X'~-&\stackrel{f}{-}&\to ~X&\\
\end{array}$$

The assumption that $K_X$ and $K_{X'}$ are nef implies that $f$ induces
an isomorphism on the complement of certain codimension $\geq2$ analytic
subsets in $X$ and $X'$ and that any exceptional divisor $E_i$ of
$\pi:Z\to X$ is also exceptional for $\pi':Z\to X'$.
Recall, that there is a positive combination $\sum n_iE_i$, i.e.\
$n_i\in\IN_{>0}$, such that $-\sum n_i E_i$ is $\pi$-ample.

Using the two direct sum decompositions
$H^2(Z,\IR)=\pi^* H^2(X,\IR)\oplus \bigoplus\IR[E_i]$ and
$H^2(Z,\IR)={\pi'}^* H^2(X',\IR) \oplus \bigoplus\IR[E_i]$,
any class $\beta\in H^2(X,\IR)$ can be written as
$\pi^*\beta={\pi'}^*\beta'+\sum a_i [E_i]$, where
$\beta'=f^*\beta$ and $a_i\in\IR$.

Now let $\alpha':=f^*(\alpha)$, where $\alpha$ is as in the proposition.
In a first step we show that all coefficients $a_i$ in
$\pi^*\alpha={\pi'}^*\alpha'+\sum a_i [E_i]$ are non-negative.

Assume that this is not the case, i.e.\
$a_1,\ldots,a_k<0$, $a_{k+1},\ldots,a_\ell\geq0$ for some $k\geq1$.
We can assume that $-(a_1/n_1)=\max_{i=1,\ldots k}\{-(a_i/n_i)\}$.
Let $C'\subset E_1$ be a general rational curve contracted under
$\pi:E_1\to X$,
which exists as  the general fibre of $\pi:E_1\to X$ is
unirational. As $C'$ is general, one has $E_i.C'\geq0$ for all
$i>1$ and hence

$$\begin{array}{clllll}
-\sum_{i=1}^{\ell}a_i(E_i.C')&\leq&-\sum_{i=1}^{k}a_i(E_i.C')&=&
-\frac{a_1}{n_1}n_1 (E_1.C')+\sum_{i=2}^k(-\frac{a_i}{n_i})n_i(E_i.C')&\\
&\leq&(-\frac{a_1}{n_1})\sum_{i=1}^k n_i (E_i.C')&\leq&(-\frac{a_1}{n_1})\sum_{i=1}^\ell n_i (E_i.C')&\\
&<&0.&&&\\
\end{array}$$

On the other hand,
$0=\pi^*\alpha.C'=({\pi'}^*\alpha'+\sum_{i=1}^\ell a_i [E_i]).C'\geq
\sum_{i=1}^\ell a_i (E_i.C')$. Contradiction.

Interchanging the r\^ole of $\alpha$ and $\alpha'$, one proves
analogously $a_i\leq0$ for all $i$. Hence, $a_i=0$ for all $i$,
i.e.\ ${\pi'}^*\alpha'=\pi^*\alpha$.

In a second step, we show that for all exceptional divisors the two
contractions $\pi:E\to X$ and $\pi':E\to X'$ coincide.
Assume that there exists a rational curve $C'\subset E_i$ contracted
by $\pi$, such that $\pi':C'\to X'$ is finite.
Since ${\pi'}^*\alpha=\pi^*\alpha$, this yields the contradiction
$0<{\pi'}^*\alpha'.C'=\pi^*\alpha.C'=0$. Analogously, one excludes the
case that there exists a rational curve $C$ contracted by $\pi'$ 
with $\pi|_C$ finite. As the fibres of $\pi'|_{E_i}$ and
$\pi|_{E_i}$ are covered by rational curves, this shows that the
two projections do coincide.

But if all contractions $\pi|_{E_i}$ and $\pi'|_{E_i}$ coincide, the
birational map $X' - - \to X$ extends to an isomorphism.\qed

\bigskip

Note that the proposition in particular applies to the case where $X$
and $X'$ admit holomorphic symplectic structures (e.g.\ $X$ and $X'$
are compact hyperk\"ahler manifolds), as in this case $K_X$ and $K_{X'}$
are even trivial.
The same arguments also prove:

\begin{corollary}\label{Hilfslemma}--- Let $f:X' - - \to X$ be
a birational map of compact complex manifolds with nef
canonical bundles. If $\alpha\in H^2(X,\IR)$ is a class
that is positive on all (rational) curves $C\subset X$, then
$\int_{C'}f^*(\alpha)>0$ for all irreducible (rational)
curves $C'\subset X'$ that are not contained in the
exceptional  locus of $f$.
\end{corollary}

\prf Under the given assumptions the arguments above show that the
coefficients $a_i$ in ${\pi'}^*\alpha'=\pi^*\alpha+\sum a_i[E_i]$
are non-negative.
If $C'\subset X'$ is an irreducible (rational) curve not
contained in the exceptional locus of $f$,
then its strict transform $\bar {C'}\subset Z$ exists and
$\int_{C'}\alpha'=\int_{\bar{C'}}{\pi'}^*\alpha'=
\int_{\bar{C'}}(\pi^*\alpha+\sum a_iE_i)>\int_{\bar{C'}}(\sum a_iE_i)\geq0$.
\qed

\bigskip

For completeness sake we recall the following results
\cite[Prop.\ 5.1,\ Cor.\ 5.2]{Huybrechts2}:

\begin{proposition}\label{Zitat}---
Let  $X$ be a compact hyperk\"ahler manifold and let
$\alpha\in H^{1,1}(X,\IR)$ be a very general element of the positive
cone $\kc_X$. Then there exist two smooth proper
families $\kx\to S$ et $\kx'\to S$ of compact hyperk\"ahler
manifolds over an one-dimensional disk $S$
and a birational map $\tilde f:\kx'- - \to \kx$ compatible with the projections
to $S$, such that $\tilde f$ induces an isomorphism
$\kx'|_{S\setminus\{0\}}\cong\kx|_{S\setminus\{0\}}$, the special
fibre $\kx_0$ is isomorphic to $X$, and $\tilde f^*\alpha$ is a K\"ahler class
on $\kx_0'$.\qed
\end{proposition}

As $S$ is contractible, $H^2(X,\IR)\cong H^2(\kx,\IR)$ and
$\tilde f^*\alpha\in H^2(\kx',\IR)=H^2(\kx_0',\IR)$.
If $\kx\leftarrow\kz\to\kx'$ resolves the birational map $\tilde f$, then
$\tilde f^*\alpha=[\Gamma]_*(\alpha)$, where
$\Gamma=\im(\kz_0\to \kx_0\times\kx_0')$.

\begin{corollary}\label{Zitatcor}---
Let $X$ be a compact hyperk\"ahler manifold and let $\alpha\in\kc_X$ be
a very general class. Then there exists a compact hyperk\"ahler
manifold $X'$ and an effective cycle $\Gamma=Z+\sum Y_i\subset X\times X'$,
such that $X\leftarrow Z\to X'$ defines a birational map $X'- - \to X$,
$[\Gamma]_*: H^*(X)\cong H^*(X')$ is a ring isomorphism compatible
with $q_X$ and $q_{X'}$, and $[\Gamma]_*(\alpha)$ is a K\"ahler
class on $X'$.\qed
\end{corollary}

By definition, the very general classes in the positive cone $\kc_X$ are
cut out by a countable number of nowhere dense closed subsets.
In fact, it suffices to assume that $q_X(\alpha,\beta)\ne0$ for all
integral class $0\ne\beta\in H^2(X,\IZ)$.
Note that there is a slight abuse of notation here. The
component $Z$ in $\Gamma$ really occurs with multiplicity one,
whereas any other component $Y_i$ may occur several times in $\sum Y_i$.
Of course, the components $Y_i$ correspond to the exceptional divisors of the
birational correspondence $\kx\leftarrow\kz\to\kx'$.

The main result of this section is the following

\begin{theorem}\label{biratgen}---
Let $X$ and $X'$ be compact hyperk\"ahler manifolds and let
$f:X' - - \to X$ be a birational map. Then there exist smooth proper
families  ${\cal X}\to S$ and
${\cal X}'\to S$ over a one-dimensional disk $S$ with the following
properties
i) The special fibres are ${\cal X}_0\cong X$ and
${\cal X}'_0\cong X'$.
ii) There exists a birational map $\tilde f:\kx'- - \to\kx$ which is an
isomorphism over $S\setminus\{0\}$, i.e.\
$\tilde f:\kx'|_{S\setminus\{0\}}\cong\kx|_{S\setminus\{0\}}$, and
which coincides with $f$ on the special fibre, i.e.\ $\tilde f_0=f$.
\end{theorem}

This generalizes a result of \cite{Huybrechts2}, where $X$ and $X'$
were assumed to be projective. In an earlier version \cite{Huybrechts}
of this result we had moreover assumed that $X$ and $X'$ are isomorphic
in codimension two. The projective arguments in the previous proofs, e.g.\
comparing sections of line bundles on $X$ and $X'$, are here replaced
by arguments from twistor theory.

\bigskip

\prf
Let $\alpha\in H^{1,1}(X,\IR)$  be a class associated with a very general
$\alpha'\in \kc_{X'}$ such that $C.\alpha'>0$ for all
(rational) curves $C'\subset X'$, e.g.\ a very general $\alpha'\in\kk_{X'}$.
Thus, $\alpha$ is a very general class in $\kc_X$ and, hence, there exist
two morphisms $\kx\to S$, $\kx'\to S$ as in Prop.\ \ref{Zitat}.
It suffices to show that $\kx'_0\cong X'$ and that under this isomorphism
$\tilde f_0$ and $f$ coincide.

The cycle $\Gamma:=\im(\kz_0\to X\times\kx_0')$ decomposes
into $\Gamma=Z+\sum Y_i$, where the $Y_i\subset X\times \kx_0'$
correspond to the exceptional divisors $D_i$
of $\kx\stackrel{\pi}{\leftarrow}\kz\stackrel{\pi'}{\to}\kx'$
and $X\leftarrow Z\to \kx_0'$ is a birational correspondence.
If the codimension of the image of $D_i$ under $\pi':D_i\to \kx_0'$
is at least two in $\kx_0'$, then $[Y_i]_*:H^2(X)\to H^2(\kx_0')$ is
trivial. If this is the case for all $i$, then
$\beta:=[\Gamma]_*(\alpha)=[Z]_*(\alpha)$. Hence, in this case,
under the birational correspondence
$\kx_0'\leftarrow Z\to X - - \to X'$ the class
$\alpha'$ on $X'$ is mapped to the K\"ahler class $\beta$ on
$\kx_0'$. Prop.\ \ref{Fujikimod} then shows that the birational map
$\kx_0'- - \to X'$ can be extended to an isomorphism. Clearly, $\tilde f_0$
corresponds to $f$.

Thus, it suffices to show that for all exceptional
divisors $D_i$ the image $\pi'(D_i)\subset \kx_0'$ has
codimension at least two.

Let $D_1,\ldots, D_k$ be those exceptional divisors for which
$\pi'(D_i)\subset\kx_0'$ is of codimension one.
We will first show that also $\pi(D_1),\ldots,\pi(D_k)\subset X$
are of codimension one. In fact,
$\pi(D_i)\subset\kx_0$ is a divisor if and only if $\pi'(D_i)\subset\kx_0'$
is one.
In order to prove this, we use that up to a non-trivial scalar
$\pi^*\sigma|_{D_i}={\pi'}^*\sigma'|_{D_i}$, where $\sigma$ and
$\sigma'$ are non-trivial two-forms on $X$ and $\kx_0$.
Hence, for any point $z\in D_i$ the homomorphisms
$\pi^*\sigma:\kt_{D_i}(z)\to\kt_X(\pi(z))\cong
\Omega_X(\pi(z))\to\Omega_{D_i}(z)$ and
${\pi'}^*\sigma':\kt_{D_i}(z)\to\kt_{\kx_0'}(\pi'(z))
\cong\Omega_{\kx'_0}(\pi'(z))\to\Omega_{D_i}(z)$ coincide.
If, for instance, $\codim(\pi(D_i))=1$, then
$\dim(\ker(\kt_{D_i}(z)\to\kt_X(z)))=1$ for $z\in D_i$
general and, therefore, $\dim(\ker\pi^*\sigma)\leq2$.
Hence, $\dim(\ker{\pi'}^*\sigma')\leq2$. This yields that
$\codim({\pi'}(D_i))=1$ or  $\ker(\kt_{D_i}(z)\to\kt_X(\pi(z)))\subset
\ker(\kt_{D_i}(z)\to\kt_{\kx_0'}(\pi'(z)))$, but the latter
is excluded, for the general fibre of $\pi:D_i\to X$ is not
contracted by $\pi'$. (In fact, this is only true for those divisors
$D_i$ on which $\pi$ and $\pi'$ differ, but if they do not differ,
then of course $\codim(\pi(D_i))=\codim(\pi'(D_i))$.)

Next, choose general irreducible (rational) curves $C_i\subset D_i$
($i=1,\ldots,k$) contracted under $\pi':D_i\to\kx_0'$.
As the exceptional locus $Sing(f^{-1})$ of the birational map
$f^{-1}: X - - \to X'$ is of codimension at
least two and $\pi(D_i)\subset X$ is a divisor,
the image $\pi(C_i)$ will not be contained in $Sing(f^{-1})$.
Moreover, we can assume that the $C_i$'s ($i=1,\ldots,k$)
do not meet $D_j$ for $j>k$. This follows from $\codim(\pi'(D_j))\geq2$ for
$j>k$.
Now, use $H^2(\kz,\IR)={\pi}^* H^2(\kx,\IR)\oplus\bigoplus \IR[D_i]=
\pi^*\Gamma(S,\underline{H}^2(X,\IR))\oplus\bigoplus \IR[D_i]$
and $H^2(\kz,\IR)={\pi'}^* H^2(\kx',\IR)\oplus\bigoplus \IR[D_i]=
\pi^*\Gamma(S,\underline{H}^2(\kx_0',\IR))\oplus\bigoplus \IR[D_i]$
and write $\pi^*\tilde\alpha={\pi'}^*\tilde\beta+\sum_{i=1}^{\ell} a_i [D_i]$,
where $\tilde\alpha\in \Gamma (S,\underline{H}^2(X,\IR))$ and
$\tilde\beta\in \Gamma (S,\underline{H}^2(\kx_0',\IR))$ are the constant
sections $\alpha$ and $\beta:=[\Gamma]_*(\alpha)$.
As in the proof of Prop.\ \ref{Fujikimod} one shows that $\beta$ K\"ahler
implies $a_i\geq0$. Also note that there exists a positive linear combination
$\sum_{i=1}^{\ell}m_iD_i$ which  is $\pi'$-negative. We may assume
that $(a_1/m_1)=\max_{i=1,\ldots,k}\{(a_i/m_i)\}$. Using Cor.\ \ref{Hilfslemma}
all this yields the contradiction

$$\begin{array}{rclcll}
0\leq\pi^*\alpha.C_1&=&({\pi'}^*\beta+\sum_{i=1}^{\ell}a_i D_i).C_1&=&\sum_{i=1}^{\ell} a_i (D_i.C_1)&\\
&=&\sum_{i=1}^{k} a_i (D_i.C_1)&\leq&\frac{a_1}{m_1}m_1(D_1.C_1)+\frac{a_1}{m_1}\sum_{i=2}^k m_i(D_i.C_1)&\\
&=&\frac{a_1}{m_1}\sum_{i=1}^k m_i(D_i.C_1)&=&\frac{a_1}{m_1}(\sum_{i=1}^\ell m_i D_i).C_1&\\
&<&0&&&\\
\end{array}$$

Hence, the varieties $\pi'(D_i)\subset \kx_0'$
are all of codimension at least two.\qed

\bigskip

Note that there are two kinds of rational curves in the proof. Those, that
sweep out a divisor in $\kx_0'$ and those that a priori do not. The former
correspond to the general fibres of the projection $D_i\to X$ of an
exceptional divisor, such that $\pi'(D_i)\subset \kx_0'$ is a divisor, and the
latter ones to those that are created by the birational correspondence
$\kx_0 - - \to X'$ itself.
In Sect.\ \ref{Kaehlercone} we will see that if $\alpha$
is positive on all rational curves, then it is automatically a K\"ahler class.
The same idea proves that the positivity of $\alpha$ on all rational curves
that sweep out a divisor is sufficient to conclude that $\alpha$
is in the ``birational K\"ahler cone''. In Sect.\ \ref{BirKcone} this is phrased
by assuming $\alpha$ to be positive on all uniruled divisors.

As mentioned in the introduction the theorem immediately
yields

\begin{corollary} --- Let $X$ and $X'$ be birational compact
hyperk\"ahler manifolds. Then the Hodge structures of $X$ and $X'$ are
isomorphic. In particular, Hodge and Betti numbers of $X$ and $X'$
coincide. Moreover, $H^*(X,\IZ)\cong H^*(X',\IZ)$ as graded rings.\qed
\end{corollary}

\begin{corollary} --- Let $X$ and $X'$ be birational compact hyperk\"ahler
manifolds of dimension $2n$. Then there exists a cycle
$\Gamma=Z+\sum Y_i\subset X\times X'$ of dimension $2n$, such that\\
i) The correspondence $X'\leftarrow Z\to X$ is the given birational
map $X'- - \to X$.\\
ii) The correspondence $[\Gamma]_*:H^*(X,\IZ)\cong H^*(X',\IZ)$
is a ring isomorphism and $[\Gamma]_*$ and $[Z]_*$ coincide on 
${H^2(X,\IZ)}$.\\
iii) The isomorphism $[\Gamma]_*$ is compatible
with $q$, i.e.\ $q_{X'}([\Gamma]_*(\, \, ))=q_X(\, \, )$.\\
iv) The isomorphism $[\Gamma]_*$ is an involution, i.e.\
$[\Gamma]_*[\Gamma]_*=\id$.\qed
\end{corollary}


\section{The K\"ahler cone}\label{Kaehlercone}

Let us recall the precise definition of the K\"ahler cone
of a compact K\"ahler manifold.

\begin{definition}\label{defkkx} --- For a complex manifold $X$ the K\"ahler
cone $\kk_X$ is the set of all classes $\alpha\in H^{1,1}(X,\IR)$
that can be represented by a closed positive $(1,1)$-form.
\end{definition}

Obviously, any ample line bundle $L$ on $X$ defines a class
$c_1(L)$ in $\kk_X$. Thus, the ample cone, i.e.\ the cone spanned
by all classes $c_1(L)$ with $L$ ample, is contained in $\kk_X$.
On the other hand, if $X$ is projective of dimension $m$ and
the canonical bundle $K_X$ is trivial, then a line bundle on $X$
is ample if and only if $\int_Xc_1^m(L)>0$ and $\int_Cc_1(L)>0$ for
any curve $C\subset X$ (cf.\ \cite[Prop.\ 6.3]{Huybrechts2}). In other words,
the ample cone can be completely described by these two conditions. For the
case of a compact hyperk\"ahler manifold (not necessarily projective)
there is a similar description for the K\"ahler cone $\kk_X$.

\begin{proposition}\label{Kcone}---
Let $X$ be a compact hyperk\"ahler manifold. A class
$\alpha\in H^{1,1}(X,\IR)$ is contained in the closure $\overline\kk_X$
of the K\"ahler cone $\kk_X\subset H^{1,1}(X,\IR)$ if and only if
$\alpha\in\overline\kc_X$ and $\int_C\alpha\geq0$ for any rational
curve $C\subset X$.
\end{proposition}

\prf Let $\alpha\in\overline\kk_X$. Then, obviously,
$\int_C\alpha\geq0$ for any curve $C\subset X$ and, since $\kk_X\subset\kc_X$,
also $\alpha\in\overline\kc_X$.

Now let $\alpha\in\overline\kc_X$ and $\int_C\alpha\geq0$ for any rational
curve $C\subset X$. If $\omega$ is any K\"ahler class,
then $\alpha+\varepsilon\omega\in\kc_X$ and
$\int_C(\alpha+\varepsilon\omega)>0$. For $\omega$ and $\varepsilon$
general, $\alpha+\varepsilon\omega$ is general as well
and if $\alpha+\varepsilon\omega\in\kk_X$
for $\omega$ general and $\varepsilon$ small, then $\alpha\in\overline\kk_X$.
Thus, it suffices to show that any general class $\alpha\in\kc_X$ with
$\int_C\alpha>0$ for any rational curve $C\subset X$ is actually K\"ahler,
i.e.\ contained in $\kk_X$.

The proof of \ref{biratgen} shows that for such a class $\alpha$ there
exists a birational map $X - - \to \kx_0'$, such that $\alpha$
corresponds to a K\"ahler class on $\kx_0'$. By Prop.\ \ref{Fujikimod}
this readily implies $X\cong\kx_0'$ and, hence $\alpha\in \kk_X$. (The
situation corresponds to the trivial case $X=X'$ in the proof of
\ref{biratgen}.).\qed

\bigskip

A weaker form of the following
result was already proved in \cite[Cor.\ 5.6]{Huybrechts2}.

\begin{corollary}--- Let $X$ be a compact hyperk\"ahler manifold
not containing any rational curve. Then $\kc_X=\kk_X$.\qed
\end{corollary}

The proposition can be made more precise: If
$\alpha\in\partial\overline\kk_X$, then $\alpha\in\partial\overline\kc_X$,
i.e.\ $q_X(\alpha)=0$, or there exists a class $0\ne\beta\in H^{1,1}(X,\IZ)$
othogonal to $\alpha$. This follows from the remark
after Cor.\ \ref{Zitatcor}.
If neither $q_X(\alpha)=0$ nor $q_X(\alpha,\beta)=0$ for some
$0\ne\beta\in H^{1,1}(X,\IZ)$, then $\alpha$ is already general and the above
proof applies directly to $\alpha$ (without adding $\varepsilon\omega$) and
shows $\alpha\in\kk_X$.
We therefore expect an affirmative answer to the following

\bigskip

\noindent
{\bf Question }--- {\it Is the K\"ahler cone $\kk_X$ of a compact hyperk\"ahler
manifold $X$ the set of classes $\alpha\in\kc_X$ with $\int_C\alpha>0$
for all rational curves $C\subset X$?}

\bigskip

In dimension two this can be proved:

\begin{corollary}\label{K3Kahler}---
Let $X$ be a K3 surface. Then the K\"ahler cone $\kk_X$ is the set
of all classes $\alpha\in\kc_X$ such that $\int_C\alpha>0$ for
all smooth rational curves $C\subset X$.
\end{corollary}

\prf
Let $C$ be an irreducible rational curve. Then $C^2\geq0$ or $C$
is smooth and $C^2=-2$. Let
$\alpha\in\kc_X$. If $C^2\geq0$, then $\alpha.C\geq0$. Hence,
the proposition shows $\overline\kk_X=\{\alpha\in\kc_X|\alpha.C\geq0
{\rm ~for~every~smooth~rational~curve}\}$. It suffices to show
that for every class $\alpha\in\partial\overline\kk_X\cap\kc_X$ there
exists a smooth rational curve $C$ with $\alpha.C=0$.
Let $\{C_i\}$ be a series of smooth irreducible rational curves,
such that $\alpha.C_i\to 0$. Let $\alpha_1=\alpha,\alpha_2.\ldots,
\alpha_{20}$ be an orthogonal base of  $H^{1,1}(X,\IR)$. We can even
require $\alpha_{i>1}^2=-1$ and $\alpha_i.\alpha_{j\ne i}=0$. For the
coefficients of $[C_i]=\lambda_i\alpha+\sum\mu_{ij}\alpha_j$ one then
concludes $\lambda_i\to 0$ and $\sum\mu_{ij}^2\to 2$. Thus,
the set of classes $[C_i]$ is contained in a compact ball and,
hence, there is only a finite number of them. The
assertion follows directly form this.\qed

\bigskip

Of course, the result is well-known, but the original arguments
are based on the Global Torelli Theorem, which seems, at least at
the moment, out of reach in higher dimensions.

The above proof shows that it would suffice to find a lower bound for
$q_X([C_i])$ for all rational curves $C_i$ such that $\int_{C_i}\alpha\to 0$.
This should also shed some light on the question whether $\partial\kk_X$ is,
away from $\partial\kc_X$, locally a finite polyhedron.

Note that also for the two other classes of compact K\"ahler manifolds
with trivial canonical bundle, compact complex tori
and Calabi-Yau manifolds, there exist descriptions of the K\"ahler
cone similar to the one for hyperk\"ahler manifolds presented above.
E.g.\ for a Calabi-Yau manifold $X$ the closure of the
K\"ahler cone $\overline\kk_X$ is the set of all
classes $\alpha\in H^2(X,\IR)$ such that $\int_X\alpha^m\geq0$ and
$\int_C\alpha\geq0$ for all (not necessarily rational) curves $C\subset X$.
This is due to the fact that on a Calabi-Yau manifold
the ample classes span the K\"ahler cone.


\section{The birational K\"ahler cone}\label{biratKaehlercone}

Let $X$ be a compact hyperk\"ahler manifold and let $f: X - - \to X'$
be a birational map to another compact hyperk\"ahler manifold $X'$.
Via the natural isomorphism $H^{1,1}(X',\IR)\cong H^{1,1}(X,\IR)$
the K\"ahler cone $\kk_{X'}$ of $X'$ can also be considered as an open
subset of $H^{1,1}(X,\IR)$ and, due to \cite[Lemma 2.6]{Huybrechts2}, in fact
as an open subset of $\kc_X$. The union of all those open subsets is
called the birational K\"ahler cone of $X$.

\begin{definition}---
The birational K\"ahler cone $\kb\kk_X$is the union
$\bigcup_{f:X - - \to X'}f^*(\kk_{X'})$, where $f:X - - \to X'$ runs
through all birational maps $X - - \to X'$ from $X$ to another compact
hyperk\"ahler manifold $X'$.
\end{definition}

Note that $\kb\kk_X$ is in fact a disjoint union of the $f^*(\kk_{X'})$,
i.e.\ if $f_1^*(\kk_{X_1})$ and $f_2^*(\kk_{X_2})$ have a non-empty
intersection then they are equal, and that $\kb\kk_X$   is in
general not a cone in $\kc_X$. As it will turn out, its closure
$\overline{\kb\kk}_X$ is a cone. Notice that this is far from being
obvious. E.g.\ if $X$ and $X'$ are birational compact hyperk\"ahler
manifolds and $\alpha$ and $\alpha'$ are
K\"ahler classes on $X$ and $X'$, respectively, there is a priori no reason
why every generic class $\alpha''$ contained in the segment joining
$\alpha$ and $\alpha'$ should be a K\"ahler class on yet another
birational compact hyperk\"ahler manifold $X''$. In this point the
theory resembles the theory of Calabi-Yau threefolds \cite{Kawamata2}.
The aim of this section is
to provide two descriptions of $\kb\kk_X$ (or rather its closure) similar
to the two descriptions of $\kk_X$ given by Def.\ \ref{defkkx} and Prop.\
\ref{Kcone}.

We begin with a geometric description of the birational K\"ahler cone
analogous to \ref{Kcone}. Here the rational curves are replaced by
uniruled divisors.

\begin{proposition}\label{BirKcone}---
Let $X$ be a compact hyperk\"ahler manifold. Then $\alpha\in H^{1,1}(X,\IR)$
is in the closure $\overline{\kb\kk}_X$ of the birational K\"ahler cone
$\kb\kk_X$ if and only if $\alpha\in\overline\kc_X$ and
$q_X(\alpha,[D])\geq0$ for all uniruled divisors $D\subset X$.
\end{proposition}

\prf
Let $\alpha\in \kb\kk_X$ and $f: X'- - \to X$ be a birational map
such that $f^*(\alpha)\in\kk_{X'}$. By \cite[Lemma 2.6]{Huybrechts2}
for any divisor $D\subset X$ we have $q_X(\alpha,[D])=
q_{X'}(f^*\alpha,f^*[D])$. Since $f^*\alpha\in\kk_{X'}$ and $f^*[D]=[f^*D]$
is effective, $q_X(\alpha,[D])>0$.
Hence, any $\alpha\in \overline{\kb\kk}_X$ is contained in
$\overline\kc_X$ and is non-negative on any effective divisor $D\subset X$ and,
in particular, on uniruled ones.

For the other direction, let $\alpha\in\kc_X$ be a general class with
$q_X(\alpha,[D])>0$ for all uniruled divisors $D\subset X$. We will show that
$\alpha\in \kb\kk_X$. This would be enough to prove the assertion.

Using Prop.\ \ref{Zitat} we find an effective cycle
$\Gamma=Z+\sum Y_i\subset X\times X'$, where $X'$ is another compact
hyperk\"ahler manifold, $Z$ defines a birational map $f: X'- - \to X$, and
$[\Gamma]_*(\alpha)\in\kk_{X'}$. It suffices to show that $[\Gamma]_*(\alpha)=
[Z]_*(\alpha)=f^*\alpha$. Clearly, if $\pi(Y_i)\subset X$ is of codimension
at least two, then $[Y_i]_*(\alpha)=0$. Here, $\pi:Y_i\to X$ and
$\pi':Y_i\to X'$ denote the first and second projection,
respectively.

Assume that $Y_1,\ldots ,Y_k$ are those components with
$\pi(Y_i)\subset X$ of codimension one
(or, equivalently $\pi'(Y_i)\subset X'$ of codimension one, as was explained
in the proof of \ref{biratgen}) and $k>0$. Firstly, we claim that 
$\sum_{i=1}^{k}[Y_i]_*(\sum_{i=1}^k[\pi(Y_i)])=
-\sum_{i=1}^kb_i'[\pi'(Y_i)]$ with $b_i'>0$.
In order to see this, recall that $\Gamma=Z+\sum Y_i$ is the
special fibre $\kz_0$ of $\kz\to S$. Hence, $\ko_\kz(\Gamma)|_\Gamma$ is
trivial and, therefore, $\ko_\kz(-\sum Y_i)|_{C_j}\cong\ko_\kz(Z)|_{C_j}$,
where $C_j$ is a generic fibre of $\pi':Y_j\to X'$. As $Y_j$ meets,
but is not contained
in $Z$, this yields $\deg(\ko_\kz(Z)|_{C_j})>0$ and hence
$-\sum_{i=1}^k(Y_i.C_j)>0$. But now
$\sum_{j=1}^k[Y_j]_*(\sum_{i=1}^k[\pi(Y_i)])=\sum_{j=1}^k(\sum_{i=1}^k(Y_i.C_j))[\pi'(Y_j)]$ and we define $b'_j:=-\sum_{i=1}^k(Y_i.C_j)>0$.
Analogously, one has $\sum_{i=1}^{k}[Y_i]_*(\sum_{i=1}^k[\pi'(Y_i)])=
-\sum_{i=1}^kb_i[\pi(Y_i)]$ with $b_i>0$.
Then, the assumption on $\alpha$ and the property
$[\Gamma]_*(\alpha)\in\kk_{X'}$ yield
the following contradiction

$$\begin{array}{l}
q_X(\alpha,\sum_{i=1}^k[\pi(Y_i)])=q_{X'}([\Gamma]_*(\alpha),[\Gamma]_*(
\sum_{i=1}^k[\pi(Y_i)]))\\
=q_{X'}([\Gamma]_*(\alpha), [Z]_*(\sum_{i=1}^k[\pi(Y_i)])+
(\sum_{i=1}^k[Y_i]_*)(\sum_{i=1}^k[\pi(Y_i)])\\
=q_{X'}([\Gamma]_*(\alpha),\sum_{i=1}^k[\pi'(Y_i)])-
q_{X'}([\Gamma]_*(\alpha),\sum_{i=1}^k b_i'[\pi'(Y_i)])\\
< q_{X'}([\Gamma]_*(\alpha), \sum_{i=1}^k[\pi'(Y_i)])\\
=q_X([\Gamma]_*[\Gamma]_*(\alpha),[\Gamma]_*(\sum_{i=1}^k[\pi'(Y_i)]))\\
=q_X(\alpha,\sum_{i=1}^k[\pi(Y_i)])+q_X(\alpha,(\sum_{i=1}^k[Y_i]_*)
(\sum_{i=1}^k[\pi'(Y_i)]))\\
< q_X(\alpha,\sum_{i=1}^k [\pi(Y_i)]).\\
\end{array}$$

Here we use that $\pi(Y_i)$ is uniruled, as it
is a component of the exceptional locus of the birational
map ${\cal X}'\leftarrow{\cal Z}\to{\cal X}$.
Indeed, any exceptional divisor $D_i$ of $\pi:{\cal Z}\to{\cal X}$
is also exceptional for $\pi':{\cal Z}\to{\cal X}'$. If the
two projections on $D_i$ differ then the images under $\pi$ and $\pi'$
are uniruled. If the two projections on $D_i$ are generically equal,
than the image $\pi(D_i)$ is contained in the uniruled image of
an exceptional divisor for which this is not the case. But in our
case the dimension of $\pi(Y_i)$, $i=1,\ldots,k$ is already maximal.
\qed

\begin{corollary}---
If $X$ is a compact hyperk\"ahler manifold not containing any uniruled
divisor, then $\kb\kk_X$ is dense in $\kc_X$.\qed
\end{corollary}

In order to obtain a description of the birational K\"ahler cone analogous
to the definition(!) of $\kk_X$ (cf. \ref{defkkx}), we need the following

\begin{lemma}\label{poscurr}---
Let $X$ be a compact K\"ahler manifold of dimension $N$.
A class $\alpha\in H^{1,1}(X,\IR)$ can be represented by a closed
positive $(1,1)$-current if and only if $\int\alpha\gamma\geq0$
for all classes $\gamma\in H^{N-1,N-1}(X,\IR)$ that can be represented
by a closed $2N-2$-form whose $(N-1,N-1)$-part is positive.
\end{lemma}

\prf
A $(1,1)$-current is a continous linear function
$T:\ka^{N-1,N-1}(X)_\IR\to \IR$. It is called
positive if $T(\gamma)\geq0$ for any positive form
$\gamma\in\ka^{N-1,N-1}(X)_\IR$.
Recall that a $(N-1,N-1)$-form is positive if and only if
it locally is a positive linear combination of forms of the type
$i\alpha_1\wedge\overline\alpha_1\wedge\ldots\wedge i\alpha_{N-1}\wedge
\overline\alpha_{N-1}$, where the $\alpha_i$'s are $(1,0)$-forms.
Note that in Def.\ \ref{defkkx} we used `positive' in the sense
of  `strictly positive'. A $(1,1)$-current
$T:\ka^{N-1,N-1}(X)_\IR\to \IR$ is closed if its extension by zero
to a linear map $\ka^{2N-2}(X)_\IR\to\IR$ is trivial on
all exact forms. Note that a priori it does not suffice to show
that the $(1,1)$-current is trivial on all exact $(N-1,N-1)$-forms.

The idea is to construct first a continous linear
function $T_0$ on all $(N-1,N-1)$-forms
$\gamma^{N-1,N-1}$, where $\gamma^{N-1,N-1}$ is the
$(N-1,N-1)$-part of a closed
$(2N-2)$-form $\gamma$. This $T_0$ shall have the properties:
{\it i)} $T_0(\gamma^{N-1,N-1})=0$ for
any exact $(2N-2)$-form $\gamma$ and {\it ii)}
$T_0(\gamma^{N-1,N-1})\geq0$ for any $\gamma^{N-1,N-1}$
that is in addition positive.
This $T_0$ will then be extended to a closed positive $(1,1)$- current.

The extension argument uses the following general result
\cite[Ch.\ II]{Bourbaki}: Let $A$ be a topological vector space,
$B\subset A$ a subspace and $C\subset A$ a convexe cone, such that
$C\cap(-C)=0$ and $B\cap\stackrel{{\scriptscriptstyle o}}{C}
\ne\emptyset$. Then, any continous
linear function $T_0:B\to \IR$ with $T_0|_{B\cap C}\geq0$
can be extended to a continous linear function $T:A\to \IR$ with
$T|_C\geq0$.

In our situation we let
$A:=\ka^{N-1,N-1}(X)_\IR/\{(d\eta)^{N-1,N-1} ~|~\eta\in\ka^{2N-3}(X)_\IR\}$
and $C$ be the image of the convexe cone of all positive
forms in $\ka^{N-1,N-1}(X)_\IR$. Clearly, $C$ is again a convex cone.
In order to see that $C\cap (-C)=0$, one argues that
for $\gamma_1,\gamma_2\in \ka^{N-1,N-1}(X)_\IR$ positive
with $\overline\gamma_1=-\overline\gamma_2$ in $A$, one has
$0\leq\int\gamma_1\omega=\int((d\eta)^{N-1,N-1}-\gamma_2)\omega=\int(d\eta-\gamma_2)\omega=-\int\gamma_2\omega\leq0$
for some $\eta\in \ka^{2N-3}(X)_\IR$,
where $\omega$ is a K\"ahler form.
Thus, $\int\gamma_1\omega=\int\gamma_2\omega=0$ and, since $\omega$ is an inner
point of the cone of all positive forms in $\ka^{1,1}(X)_\IR$, this implies
$\gamma_1=\gamma_2=0$. Last but not least, if
$B:=\{\gamma^{N-1,N-1}\in A~|~d(\gamma)=0\}$ then $B\cap C$ contains
$\omega^{N-1}$, which is an inner point.

Therefore, in order to construct the closed positive
current representing a class $\alpha\in\kc_X$ one has to define $T_0:B\to \IR$
such that $T_0$ is non-negative on $B\cap C$. The obvious choice
is $\gamma^{N-1,N-1}\mapsto\int\gamma^{N-1,N-1}\alpha=\int\gamma\alpha$,
which only depends on the cohomology class of $\alpha$
as long as $\gamma$ is closed and $\alpha$ is a closed
$(1,1)$-form representing its cohomology class. It is well
defined on $B$, because $\int\alpha\gamma^{N-1,N-1}=0$ for all exact $\gamma$.
It remains to show that $\int\alpha\gamma\geq0$ for
any positive $\gamma^{N-1,N-1}\in\ka^{N-1,N-1}(X)_\IR$ with
$d(\gamma)=0$.
But this is ensured by the assumption.\qed

\bigskip

In analogy to the definition of $\kk_X$ we can now prove the

\begin{proposition}\label{descrbir}---
Let $X$ be a compact hyperk\"ahler manifold.
The closure $\overline{\kb\kk}_X$ of the birational K\"ahler cone $\kb\kk_X$
is the set of all $\alpha\in H^{1,1}(X,\IR)$, such that
$\alpha(\sigma\bar\sigma)^{n-1}\in H^{2n-1,2n-1}(X,\IR)$ can be represented
by a closed $(4n-2)$-form whose $(2n-1,2n-1)$-part is positive.
\end{proposition}

\prf Assume that $\alpha\in\overline{\kb\kk}_X$, but that
$\alpha(\sigma\bar\sigma)^{n-1}\in H^{2n-1,2n-1}(X,\IR)$ cannot be represented
by a closed $(4n-2)$-form with
positive $(2n-1,2n-1)$-part. Since the set of classes
$\gamma\in H^{2n-1,2n-1}(X,\IR)$ that can be represented by such forms
is convex, there exists a class $\beta\in H^{1,1}(X,\IR)$, which is
positive on those but negative on $\alpha(\sigma\bar\sigma)^{n-1}$.
Due to the above lemma $\beta$ can be represented by a closed
positive $(1,1)$-current.

Let $f:X' - - \to X$ be a birational map from another compact hyperk\"ahler
manifold $X'$. Then also $f^*\beta$ can be represented by a closed
positive $(1,1)$-current (cf. \cite{Fujiki3}) and, therefore,
$f^*\beta$ is positive on $\kk_{X'}$ with respect to $q_{X'}$.
Hence, $\beta$ must be non-negative on $\overline{\kb\kk}_X$, which
contradicts $q_X(\alpha,\beta)=\int\alpha\beta(\sigma\bar\sigma)^{n-1}<0$.

On the other hand, if $\alpha(\sigma\bar\sigma)^{n-1}$ can be represented
by a closed $(4n-2)$-form with positive $(2n-1,2n-1)$-part, then
$\int_D\alpha(\sigma\bar\sigma)^{n-1}\geq0$ for all divisors $D\subset X$.
Prop.\ \ref{BirKcone} then in particular yields $\alpha\in\overline{\kb\kk}_X$.
\qed

\bigskip

\begin{remark}-- It might be that in Lemma \ref{poscurr} it suffices to
test the class $\alpha$ on  those $\gamma$ that can be represented by closed
positive $(N-1,N-1)$-forms. The result in Prop.\ \ref{descrbir} would
change accordingly.
\end{remark}

\begin{corollary}---
The closure of the birational K\"ahler cone $\overline{\kb\kk}_X$
is dual to the cone of classes $\alpha\in H^{1,1}(X,\IR)$ that can be
represented by closed positive $(1,1)$-currents.\qed
\end{corollary}

\section{Curves of negative square}\label{negcurve}

For K3 surfaces it is well-known that an irreducible curve with negative
self-intersection is a smooth rational curve. As the
quadratic form $q_X$ on an arbitrary higher-dimensional
compact hyperk\"ahler manifold is replacing the intersection
pairing, one may wonder whether curves (or divisors) of negative square also
have special geometric properties in higher dimensions.

First note that the quadratic form $q_X$ on $H^{1,1}(X)$ for
a compact hyperk\"ahler manifold $X$ also induces a quadratic
form on $H^{2n-1,2n-1}(X)$, also denoted by $q_X$.
Here, one can either use the natural duality isomorphism
$H^{1,1}(X)\cong H^{2n-1,2n-1}(X)\dual$ or the isomorphism
$L:H^{1,1}(X)\cong H^{2n-1,2n-1}(X)$ given by taking the product with
$(\sigma\bar\sigma)^{n-1}$, where $0\ne\sigma\in H^0(X,\Omega_X^2)$.
Thus, we can speak of curves $C\subset X$ with square $q_X([C])$.
For curves of negative square Prop.\ \ref{Kcone} yields:

\begin{proposition}\label{prop51}---
Let $X$ be a compact hyperk\"ahler manifold and
$C\subset X$ be a curve with
$q_X([C])<0$. Then there exists a class $\alpha\in\kc_X$ and an irreducible
rational curve $C'\subset X$ such that $\alpha.C<0$ and $\alpha.C'<0$.
\end{proposition}

\prf
Let $\beta:=L^{-1}([C])$. Then, $q_X(\beta)=q_X([C])<0$.
If $\gamma\in\kk_X$, then $q_X(\gamma,\beta)=\int_C\gamma>0$. Therefore,
the class $\alpha:=(\gamma-(2q_X(\gamma,\beta)/q_X(\beta))\beta)$
satisfies $\int_C\alpha=q_X(\alpha,\beta)=-\gamma.C<0$,
$q_X(\alpha)=q_X(\gamma)>0$, and $q_X(\alpha,\gamma)>0$.
Thus, the class $\alpha$ is contained in the positive
cone $\kc_X$, but $\alpha\notin\overline\kk_X$.
By \ref{Kcone} there exists an irreducible rational curve $C'\subset X$,
such that $\alpha.C'<0$.\qed

\bigskip

\begin{corollary}---
In particular, if there are no irreducible rational curves of negative
square in $X$, then all(!) curves have non-negative square.\qed
\end{corollary}

Also, if $H^{2n-1,2n-1}(X,\IQ)=\IQ[C]$, where $C\subset X$ is of negative
square, then there exists a rational curve $C'\subset X$
with $\IQ[C']=\IQ[C]$. This (and the deformation theory of
hyperk\"ahler manifolds) suggests an affirmative answer to the following

\bigskip

\noindent
{\bf Question} --- {\it Let $C\subset X$ be an irreducible curve in a compact
hyperk\"ahler manifold $X$ with $q_X([C])<0$. Does there exist a rational
curve $C'\subset X$ such that $\IQ[C]=\IQ[C']$ in $H^{2n-1,2n-1}(X,\IQ)$?}

\bigskip

Something slightly weaker can in fact be proved:

\begin{proposition}\label{Uni}---
Let $C$ be an irreducible curve in a compact hyperk\"ahler manifold
$X$. If $q_X([C])<0$, then $C$ is contained
in a uniruled subvariety $Y\subset X$.
\end{proposition}

\prf As in the previous proof one finds a very general class
$\alpha\in\kc_X$ with $\alpha.C<0$.
Therefore, there exist two families $\kx\to S$ and $\kx'\to S$ as
in Prop.\ \ref{Zitat}. Let $\kx\leftarrow\kz\to\kx'$
be a birational correspondence induced by the isomorphism
$\kx|_{S\setminus\{0\}}\cong \kx'|_{S\setminus\{0\}}$.
As the exceptional locus of the birational map ${\cal X}- - \to {\cal X}'$
is uniruled ($K_{\cal X}$ and $K_{{\cal X}'}$ are trivial),
it suffices to show that $C\subset X=\kx_0$ is in the exceptional
locus of $\kx\leftarrow \kz\to \kx'$. If not, the strict transform
$\bar C\subset\kz$ of $C$ exists and
$\tilde\alpha.C=\pi^*\tilde\alpha.\bar C=
({\pi'}^*\tilde\alpha'+\sum a_i [D_i]).\bar C$ (cf. Cor.\ \ref{Hilfslemma}).
Here, we use the notation of the proof of Thm.\ \ref{biratgen}, i.e.\ the
$D_i$'s are the exceptional divisors of $\kz\to\kx$,
$\tilde\alpha'_0=\alpha'$ is a K\"ahler class on $\kx_0'$, and
the coefficients $a_i$'s are positive.
This yields the contradiction:
$0>\alpha.C={\pi'}^*\tilde\alpha'.\tilde C+\sum a_i (D_i.C)>0$.\qed

\bigskip

If one could replace  uniruled by unirational in Prop.\ \ref{Uni},
we would obtain  $\IQ[C]=\IQ[\sum a_iC_i]$, where the
$C_i$'s are rational curves, but the $a_i\in\IZ$ might be negative
(cf.\ \cite{Kollar}).

Analogously to Prop.\ \ref{prop51} one can use the description of the
birational K\"ahler cone (Prop.\ \ref{BirKcone}) to prove:

\begin{proposition}---
Let $X$ be a compact hyperk\"ahler manifold and let $D\subset X$ be
an effective divisor such that $q_X([D])<0$. Then there exists a class
$\alpha\in\kc_X$ and a uniruled divisor $D'\subset X$ such that
$q_X(\alpha,[D])<0$ and $q_X(\alpha,[D'])<0$. In particular, if $X$
does not contain any uniruled divisor then there does not exist any divisor
$D\subset X$ with $q_X([D])<0$.\qed
\end{proposition}


\noindent
{\small Mathematisches Institut\\
Universit\"at zu K\"oln\\
Weyertal 86-90\\
50931 K\"oln, Germany\\}

\end{document}